\documentclass[a4paper, 10pt, reqno]{amsart}

\usepackage{ amssymb, amsmath, amsthm}
\usepackage{graphicx, psfrag, setspace, subfigure, gensymb}
\usepackage{amsfonts}
\usepackage{bbm}
\usepackage{fix-cm}
\usepackage{comment}
\usepackage{tikz,tikz-cd}
\usetikzlibrary{shapes}
\usetikzlibrary{matrix,arrows,decorations.pathmorphing,decorations.pathreplacing}
\tikzset{commutative diagrams/diagrams={baseline=-2.5pt},commutative diagrams/arrow style=tikz}
\usepackage[colorlinks]{hyperref}
\usepackage{float}
\usepackage{setspace} % needed for redefined \marginpar

\newcommand\Z{\mathbb Z}
\newcommand\C{\mathbb C}

\newcommand\N{\mathbb N}

\newcommand{\cO}{\mathcal{O}}

\newcommand{\cR}{\mathcal{R}}

\newcommand\into{\hookrightarrow}
\newcommand\To{\longrightarrow}

\newcommand\Mat{\operatorname{Mat}}

\renewcommand\P{\mathbb P}
\newcommand\Gr{\operatorname{Gr}}

\newcommand{\Wedge}{\mbox{\scalebox{1.2}{$\wedge$}}}
\newcommand{\mat}[1]{\begin{pmatrix}#1\end{pmatrix}}

\newcommand{\al}[1]{\begin{align*}#1\end{align*}}

\newcommand{\beq}[1]{\begin{equation}\label{#1} }
\newcommand{\eeq}{\end{equation}}
\newcommand{\pgap}{\vspace{5pt}}

\theoremstyle{plain}

\newtheorem{thm}[equation]{Theorem}

\theoremstyle{remark}
\newtheorem{rem}[equation]{Remark}

\theoremstyle{definition}

\newtheorem{eg}[equation]{Example}

\newtheorem{claim}[equation]{Claim}

\makeatletter \@addtoreset{equation}{section} \makeatother

\let\oldtocsection=\tocsection
\let\oldtocsubsection=\tocsubsection
\let\oldtocsubsubsection=\tocsubsubsection
\renewcommand{\tocsection}[3]{\hspace{0em}\oldtocsection{#1}{#2}{#3}}
\renewcommand{\tocsubsection}[3]{ \hspace{1em} \oldtocsubsection{#1}{\small{#2}}{\small{#3}} }
\renewcommand{\tocsubsubsection}[3]{\hspace{2em}\oldtocsubsubsection{#1}{\small{#2}}{\small{#3}}}

\newcommand{\marginparstretch}{0.6}
\let\oldmarginpar\marginpar
\renewcommand\marginpar[1]{\-\oldmarginpar[\framebox{\setstretch{\marginparstretch}\begin{minipage}{\marginparwidth}{\raggedleft\scriptsize #1}\end{minipage}}]{\framebox{\setstretch{\marginparstretch}\begin{minipage}{\marginparwidth}{\raggedright\scriptsize #1}\end{minipage}}}}

\newtheorem{examplex}[equation]{Example}
\newenvironment{egu}
  {\pushQED{\qed}\examplex}
  {\popQED\endexamplex}

\newcommand{\aand}{\quad\quad\mbox{and}\quad\quad}
\newcommand\ie{\emph{i.e.}~}

\begin{document}

\title{A short guide to GKZ}
\author{Ed Segal}

\maketitle

\begin{abstract} These notes are a brief summary of the main results from the book `Discriminants, Resultants and Multidimensional Determinants' by Gelfand--Kapranov--Zelevinsky. We sketch the key ideas involved in the proofs, using as little technical background as possible. 
\end{abstract}

\tableofcontents

\section{Introduction}\label{sec.intro}

The purpose of these notes is to provide a short ``user's guide'' to the wonderful book \emph{Discriminants, Resultants and Multidimensional Determinants} by Gelfand, Kapranov and Zelevinsky \cite{GKZ}. This book is full of interesting results and examples, very clearly explained. However, since the main results are developed in great depth over several hundred pages, I felt there might be some value in writing an outline of what I understand of their central argument, in case it helps others navigate the book themselves.

 I have tried to keep it as brief and as non-technical as possible while still providing all the essential ideas for understanding the proof. This means that many topics are skipped or covered very briefly, and of course my choice of what counts as `essential' is highly subjective. Fortunately, the reader who wants to know more about any topic already knows where they should look.
\pgap

\textbf{Acknowledgements.} I would like to thank Hiroshi Iritani, for originally introducing me to this subject; Alex Kite and Will Donovan, for encouraging me to actually read the book; and Michela Barbieri, for very helpful feedback on a draft.

\subsection{The discriminant of a cubic}\label{cubic}

We all learned in kindergarten that the discriminant of a quadratic equation
$$aX^2 + bX + c$$
is $\Delta_2=b^2-4ac$. In the 3-dimensional space of quadratics there is a hypersurface $\{\Delta_2=0\}$ where the quadratic has a repeated root. 

The analogous question for cubic equations is much harder. The discriminant of a cubic
$$aX^3+bX^2 + cX+d $$
is the polynomial
\beq{cubicdis} \Delta_3 = b^{2}c^{2}-4ac^{3}-4b^{3}d-27a^{2}d^{2}+18abcd \eeq
\ie we have $\Delta_3(a,b,c,d)=0$ iff the cubic has a repeated root. This formula was probably well-known to mathematicians of earlier centuries, but I suspect few modern mathematicians could write it down without googling. The discriminant of a quartic, $\Delta_4$, is a polynomial with 16 terms. 
\pgap

The most important (in my view) result of \cite{GKZ} is that in fact a lot of information about discriminants can be uncovered using a simple combinatorial game. Let's illustrate this procedure for the case of cubics. We start by drawing the \emph{Newton polytope} of the cubic - the convex hull of the monomials involved. This is just an interval:
$$
\begin{tikzpicture}[scale=1.5]
\draw (0,0)--(3,0);
\filldraw (0,0) circle (1pt) node[above, yshift=2pt]{$aX^3$};
\filldraw (1,0) circle (1pt)node[above, yshift=2pt]{$bX^2$};
\filldraw (2,0) circle (1pt)node[above, yshift=2pt]{$cX$};
\filldraw (3,0) circle (1pt)node[above, yshift=2pt]{$d$};
\end{tikzpicture}
$$
Next we subdivide this interval into smaller intervals. There are exactly four ways we can do this:
$$
\begin{tikzpicture}[scale=.7]
\draw (0,0)--(3,0);
\filldraw (0,0) circle (1pt) node[above, yshift=2pt]{$a$};
\filldraw (1,0) circle (1pt)node[above, yshift=2pt]{$b$};
\filldraw (2,0) circle (1pt)node[above, yshift=2pt]{$c$};
\filldraw (3,0) circle (1pt)node[above, yshift=2pt]{$d$};
\end{tikzpicture}
\quad\quad
\begin{tikzpicture}[scale=.7]
\draw (0,0)--(3,0);
\filldraw (0,0) circle (1pt) node[above, yshift=2pt]{$a$};
\draw (1,0)--(1,0.1);
\filldraw (2,0) circle (1pt)node[above, yshift=2pt]{$c$};
\filldraw (3,0) circle (1pt)node[above, yshift=2pt]{$d$};
\end{tikzpicture}
\quad\quad
\begin{tikzpicture}[scale=.7]
\draw (0,0)--(3,0);
\filldraw (0,0) circle (1pt) node[above, yshift=2pt]{$a$};
\filldraw (1,0) circle (1pt)node[above, yshift=2pt]{$b$};
%\filldraw (2,0) circle (1pt)node[above, yshift=2pt]{$c$};
\draw (2,0)--(2,0.1);
\filldraw (3,0) circle (1pt)node[above, yshift=2pt]{$d$};
\end{tikzpicture}
\quad\quad
\begin{tikzpicture}[scale=.7]
\draw (0,0)--(3,0);
\filldraw (0,0) circle (1pt) node[above, yshift=2pt]{$a$};
\draw (1,0)--(1,0.1);
\draw (2,0)--(2,0.1);
%\filldraw (1,0) circle (1pt)node[above, yshift=2pt]{$b$};
%\filldraw (2,0) circle (1pt)node[above, yshift=2pt]{$c$};
\filldraw (3,0) circle (1pt)node[above, yshift=2pt]{$d$};
\end{tikzpicture}
$$
The first subdivision is into three intervals each of length 1, the second and third use an interval of length 2 and another of length 1, and the fourth subdivision uses a single interval of length 3. 

For each of these subdivisions we can produce a monomial in the variables $a,b,c,d$. Here's the rule: for each subinterval, consider the expression
\beq{game}\big(\mbox{length of subinterval}\times \mbox{endpoints}\big)^{\mbox{\small{length}}} \eeq
and then multiply over each subinterval appearing. So the first of our subdivisions gives the monomial
$$ ab\times bc \times cd = ab^2c^2 d$$
and the second subdivision gives:
$$ (2ac)^2 \times cd= 4a^2c^3d$$
The third and fourth subdivisions give monomials $4ab^3d^2$ and $27a^3d^3$ respectively. We have nearly recovered the discrimant of a cubic! If we divide by the common factor of $ad$, and ignore the sign\footnote{The issue of signs is delicate, it is addressed in \cite{GKZ} but we will not cover it in these notes.}, we have found four of the five terms in \eqref{cubicdis}. 
\pgap

What of the fifth term in $\Delta_3$? To understand why this term is different from the others we need to draw the Newton polytope of $\Delta_3$ itself. At first sight this polytope lives in four dimensions so is hard to draw, but there are some obvious symmetries that we can use to our advantage. Namely (i) scaling by an overall constant will clearly not affect whether or not our cubic has a repeated root, and (ii) neither will scaling the variable $X$ by some factor. So the space of cubics $\C^4_{a,b,c,d}$ carries an action of a torus $(\C^*)^2$ under which $\Delta_3$ is invariant; the first factor acts with weights $(1,1,1,1)$ and the second with weights $(3,2,1,0)$. If we allow ourselves to pass to rational functions in $a,b,c,d$ then it's easy to compute that the invariants are freely generated by
$$\alpha = \frac{b^2}{ac} \aand \beta = \frac{c^2}{bd}$$
so in fact $\Delta_3$ must be (proportional to) a rational function in $\alpha, \beta$. And indeed, a little manipulation reveals that
$$\frac{\Delta_3}{a^2d^2} = \alpha^2\beta^2  - 4\alpha\beta^2 - 4 \alpha^2\beta - 27 + 18\alpha\beta$$
which has Newton polytope:
$$
\begin{tikzpicture}[scale=1.5]
\draw (0,0)--(2,1)--(2,2)--(1,2)--(0,0);
\filldraw (0,0) circle (1pt) node[above, xshift=-2pt, yshift=2pt]{$1$};
\filldraw (1,1) circle (1pt)node[above, yshift=2pt]{$\alpha\beta$};
\filldraw (2,1) circle (1pt)node[above, xshift=12pt, yshift=2pt]{$\alpha^2\beta$};
\filldraw (1,2) circle (1pt)node[above, yshift=2pt]{$\alpha\beta^2$};
\filldraw (2,2) circle (1pt)node[above, yshift=2pt]{$\alpha^2\beta^2$};
\end{tikzpicture}$$
We can see that our combinatorial game has found the terms of $\Delta_3$ that live at the vertices of the Newton polytope. We will refer to these as \emph{extremal terms}. The term $18\alpha\beta$ lives in the interior. 

\subsection{More general discriminants}\label{general}

It's easy to check that the discriminant of a quadratic appears when we play the same game starting with an interval of length 2. There are two possible subdivisions and they give us terms $ab^2c$ and $4a^2c^2$, then dividing by $ac$ (and inserting a minus sign) gives us $\Delta_2$. Here all terms are extremal since the Newton polytope of $\Delta_2$ is very small.

 If we play the game for an interval of length 4 then we begin to discover the discriminant of a quartic equation; it begins
$$\Delta_4 = b^2c^2d^2 \pm 4ac^3d^2e \pm 256a^3e^3 \pm ... $$
and one can easily write down another five terms coming from the remaining subdivisions. These eight terms will be the extremal terms of $\Delta_4$.
\pgap

More generally we can consider discrimants of polynomials in more than one variable. The generalization of `having a repeated root' is `defining a singular hypersurface'. 

\begin{egu} \label{quadsurf}
Consider the following polynomial in two variables:
$$f(X,Y) = a + bX + cY + dXY $$
It defines a hypersurface $(f)\subset \C^2$, and it's easy to calculate that $(f)$ has a singularity iff:
$$ad-bc = 0$$
So the discriminant of $f$ is $ad-bc$. This is one of the rare examples where it is easy to find the discriminant directly.

The Newton polytope of $f$ is a square. The generalization of `subdividing into intervals' is `triangulating', \ie subdividing into simplices.\footnote{We will only ever consider \emph{coherent} triangulations \cite[p218]{GKZ}. In these notes `triangulation' means `coherent triangulation'.} Since there are two possible triangulations 

%\footnote{But also partly because of some unfortunate numbering conventions in the typesetting. There are, for example, four separate results labelled Proposition 3.1.}

$$
\begin{tikzpicture}[scale=1]
\draw (0,0)--(1,0)--(1,1)--(0,1)--(0,0);
\draw (0,0)--(1,1);
\filldraw (0,0) circle (1pt) node[above, xshift=-5pt, yshift=0pt]{$a$};
\filldraw (1,0) circle (1pt)node[above, xshift=5pt, yshift=0pt]{$b$};
\filldraw (0,1) circle (1pt)node[above, xshift=-5pt,yshift=0pt]{$c$};
\filldraw (1,1) circle (1pt)node[above, xshift=5pt,yshift=0pt]{$d$};
\end{tikzpicture}
\hspace{2cm}
\begin{tikzpicture}[scale=1]
\draw (0,0)--(1,0)--(1,1)--(0,1)--(0,0);
\draw (0,1)--(1,0);
\filldraw (0,0) circle (1pt) node[above, xshift=-5pt, yshift=0pt]{$a$};
\filldraw (1,0) circle (1pt)node[above, xshift=5pt, yshift=0pt]{$b$};
\filldraw (0,1) circle (1pt)node[above, xshift=-5pt,yshift=0pt]{$c$};
\filldraw (1,1) circle (1pt)node[above, xshift=5pt,yshift=0pt]{$d$};
\end{tikzpicture}
$$
the combinatorics predicts that the disciminant of $f$ will be:
$$ (abd)(acd) \pm (abc)(bcd) = abcd(ad \pm bc) $$
Here we are using the formula \eqref{game} generalized to:
\beq{game2}(\mbox{volume of simplex}\times\mbox{vertices})^{\mbox{\small{volume}}} \eeq
Apart from the factor of $abcd$ the game works, and there are no possible interior terms. 

\end{egu}

If we have $d$ variables $X_1,..., X_d$ and we fix a finite set 
$$A \subset \N^d$$
of monomials in these variables then there is a corresponding generic polynomial $f_A$ with these terms, depending on $n=|A|$ parameters $a_1,..., a_n$. For the hypersurface $(f_A)\subset \C^d$ to be singular we need $d+1$ polynomials to vanish simultaneously, $f_A$ and $\partial_i f_A, \forall i$. So we expect that $(f_A)$ is smooth for a generic choice of parameters, and singular for a codimension 1 locus in $\C^n$. Assuming that this is true we define the \emph{discriminant} of $A$ to be the polynomial 
$$\Delta_A \in \C[a_1,..., a_n]$$
cutting out the locus where $(f_A)$ is singular. As we see already for cubics in one variable, computing this polynomial $\Delta_A$ is going to be a difficult problem.

\begin{rem}\label{resultants}
Finding discriminants is a special case of a more general problem called `resultants'. Fix $d+1$ sets $A_0,..., A_d$ of monomials, all in the same set of variables $X_1,.., X_d$. Then there is a corresponding list of polynomials $g_0,..., g_d$ and for generic values of the co-efficients we expect that the polynomials will have no common root. But there will be (we expect) a hypersurface in parameter space where the $g_i$'s \emph{do} have a common root. The polynomial cutting out this hypersurface is called the \emph{resultant}.  

Obviously if we can compute resultants then we can compute discriminants, by setting $g_0=f$ and $g_i=\partial_{X_i} f$. This idea will be important in Sections \ref{associated} and \ref{fromchow}.
\end{rem}

In some examples the naive dimension counting argument fails, and the subset of the parameter space $\C^n$ where $(f_A)$ is singular is not a hypersurface. 

\begin{eg} Let $d=1$ and $A=\{X^4, X^3, X^2\}$ so $f_A = aX^4 + bX^3 + cX^2$. Obviously $(f_A)$ is always singular at the origin, regardless of the values of $a,b,c$.  However, if we remove the origin then $(f_A)$ is singular iff $b^2-4ac=0$. 

\end{eg}
This example shows that it is better to consider $(f_A)$ as a hypersurface in the torus $(\C^*)^d$ instead of $\C^d$. If we do this then $\Delta_A$ is not affected by translations of the set $A$, and indeed we can allow $f_A$ to be a Laurent polynomial, \ie with $A\subset \Z^d$ instead of $\N^d$. 

The next example is a little more subtle.

\begin{egu}\label{eg.othercomponentsneeded}
Let $d=2$ and $A=\{1, X, X^2, Y\}$, so:
$$f_A = a + bX + cX^2 + dY$$
Obviously $(f_A)\subset (\C^*)^2$ is singular iff the parameters lie in the subset 
$$\{d= b^2-4ac=0\}$$
 which has codimension two in  $\C^4$. 

What does our combinatorial game predict for this example? The convex hull of $A$ has two possible triangulations
$$
\begin{tikzpicture}[scale=1]
\draw (0,0)--(1,0)--(2,0)--(0,1)--(0,0);
\draw (0,1)--(1,0);
\filldraw (0,0) circle (1pt) node[below, yshift=0pt]{$a$};
\filldraw (1,0) circle (1pt)node[below, yshift=0pt]{$b$};
\filldraw (2,0) circle (1pt)node[below,yshift=0pt]{$c$};
\filldraw (0,1) circle (1pt)node[above, yshift=0pt]{$d$};
\end{tikzpicture}
\hspace{2cm}
\begin{tikzpicture}[scale=1]
\draw (0,0)--(1,0)--(2,0)--(0,1)--(0,0);
\filldraw (0,0) circle (1pt) node[below, yshift=0pt]{$a$};
\draw (1,0)--(1,0.1);
\filldraw (2,0) circle (1pt)node[below,yshift=0pt]{$c$};
\filldraw (0,1) circle (1pt)node[above, yshift=0pt]{$d$};
\end{tikzpicture}
$$
so the game predicts that the discriminant of $A$ should be:
$$ (abd)(bcd) \pm (2acd)^2 = acd^2(b^2 \pm 4ac) $$

It is apparent that this expression is \emph{not} detecting the disciminant of $f_A$ itself, but rather of the polynomial
$$f_{A'} = a + bX + cX^2$$
associated to the bottom face $A'\subset A$. 

\end{egu}

In light of this we define a polynomial $E_A$, which \cite{GKZ} call the \emph{principal $A$-determinant}, by setting:
\beq{principalAdet}E_A = \prod_{\mbox{\tiny{faces }}A'\subset A} (\Delta_{A'})^{m_{A'}} \eeq
Here:
\begin{itemize}\setlength{\itemsep}{3pt}
\item The product runs over all faces of $A$ of all dimensions, including $A$ itself.
\item We formally declare $\Delta_{A'}\equiv 1$ if the locus where $(f_{A'})$ is singular is not a hypersurface.
\item The $m_{A'}$'s are certain multiplicities that we will not address in these notes.
\end{itemize}
And then their main result is:

\begin{thm}\cite[p302]{GKZ}\footnote{We're citing results in \cite{GKZ} by page number because we found their numbering conventions challenging. There are, for example, four separate results labelled Proposition 3.1.} \label{mainthm}The combinatorial game described above computes the extremal terms of $E_A$.
\end{thm}

We will sketch the proof of this theorem in Section \ref{proof}. But before we start let's make a few elementary remarks:

\begin{itemize}\setlength{\itemsep}{5pt}
\item The discriminant locus $\{\Delta_A=0\}$ is now just one of the irreducible components\footnote{Each hypersurface $\{\Delta_{A'}=0\}$ is irreducible by a simple geometric argument \cite[p15]{GKZ}.} of a larger variety $\{E_A=0\}$. We call it the \emph{principal component}, and it may be empty as in Example \ref{eg.othercomponentsneeded}. It's often sensible to think of the larger set $\{E_A=0\}$ as the correct `discriminant locus' associated to $A$.

\item Only certain faces $A'\subset A$ will contribute non-empty components of the discriminant. They can be identified by a simple combinatorial criterion, see \emph{e.g.}~\cite[Sec.~4.2.2]{KS}. 

\item If $A'$ is a vertex of $A$ then $f_{A'}$ is a single monomial, and $\Delta_{A'}$ is just the coefficient of this monomial. So each variable lying at a vertex of $A$ appears as a factor of $E_A$. For example, the principle $A$-determinant for a quadratic in one variable will be
\beq{eq.quadricEA} E_A= ac(b^2-4ac) \eeq
(all the multiplicities $m_{A'}$ are 1 in this example). This explains the additional factors that appeared in all our examples above.
\end{itemize}

\section{Proof of the main theorem}\label{proof}

\subsection{Projective duality}

Discriminants are closely related to the geometric subject of projective duality, which is almost as classical and even more subtle. 
\pgap

Let $V\subset \P^n$ be a smooth projective variety. Write $(\P^n)^\vee$ for the dual projective space of hyperplanes in $\P^n$. The \emph{projective dual} of $V$ is the subvariety:
$$V^\vee = \{H \in (\P^n)^\vee\; |\; H \mbox{ is tangent to } V \} $$ 
Note that $V^\vee$ is not intrinsic to $V$, it depends on the embedding $V\hookrightarrow\P^n$. Also we can extend the definition to singular $V$ by declaring $V^\vee$ to be the closure of the set of hyperplanes which are tangent to a smooth point of $V$. 

 The key property of projective duality is that it is indeed a duality, \ie the map  $V \mapsto V^\vee$ is an involution on the set of subvarieties of $\P^n$. This is called the \emph{biduality theorem} and requires a little work to prove \cite[p15]{GKZ}. 
\pgap

By a dimension count, as in Section \ref{general}, we see that we expect $V^\vee$ to be a hypersurface in $(\P^n)^\vee$. This is regardless of the dimension of $V$.

\begin{rem}
 At first sight this fact seems odd given the biduality theorem; we have an involution $V\leftrightarrow V^\vee$ on subvarieties of $\P^n$ under which almost every subvariety maps to a hypersurface. However, if $V$ is not a hypersurface then $V^\vee$ is a rather special variety: it must be ruled by linear subspaces (because if $\operatorname{codim} V=c$ then every smooth point of $V$ has a $\P^{c -1}$ of tangent hyperplanes). A generic deformation of $V^\vee$ will lose this property so must be projectively-dual to a hypersurface, not to a deformation of $V$. So the moduli space of $V^\vee$ typically contains the moduli space of $V$ as a subvariety of lower dimension.  See Example \ref{twistedcubic2}.
\end{rem}

Now we relate projective duality to disciminants. Fix, as in the previous section, a lattice $\Z^d$ and a finite subset $A\subset \Z^d$ of size $n$. We will explain how to construct a projective toric variety
$$V_A \subset \P^{n-1}$$
such that the discriminant $\Delta_A$ is exactly the equation of the projective dual of $V_A$.\footnote{We shall see in Section \ref{mirror} that $V_A$ is one of various different toric varieties associated to the data of $A$.}
\pgap

There a few different ways to express this construction. The simplest is to regard our set of Laurent monomials $A$ as defining an embedding
$$ (\C^*)^d \into \P^{n-1}$$
and take $V_A$ to be the closure of the image.

\begin{egu}\label{quadcurve1}
Returning to quadratics in one variable, set $A=\{X^2, X, 1\}$. This defines an embedding:
\al{ \C^* & \into \P^2 \\ x &\; \mapsto x^2\!:\! x \!:\! 1 }
The closure of the image is obviously the quadric curve $(\alpha\gamma- \beta^2)$, and the embedding extends to 
\begin{align}V_A \cong \P^1 & \into \P^2 \label{eq.quadembed} \\ \nonumber
 x\!:\! y  \;& \mapsto \; x^2\!:\! xy \!:\! y^2 \end{align}
\end{egu}

More formally, we can embed $\Z_d$ into $\Z_{d+1}$ as an affine hypersurface at height 1, and view $A$ as a subset of $\Z^{d+1}$.  Then we can consider the cone
$$\Gamma_A = \langle A\rangle_{\N} \subset \Z^{d+1} $$
\emph{i.e.} the submonoid of $\Z^{d+1}$ generated by $A$. Then the monoid ring $\C[\Gamma_A]$ is a graded ring generated in degree 1. We define $V_A$ to be the associated projective variety:
$$V_A = \operatorname{Proj} \C[\Gamma_A] $$
By construction this is a toric variety, equipped with an embedding into $\P^{n-1}$.

\begin{egu}\label{twistedcubic}
 If we set $A=\{ X^3, X^2, X, 1\}$ then it's easy to verify that 
\beq{twistedcubicideal}\C[\Gamma_A] = \C[\alpha, \beta, \gamma, \delta] \,/\,(\alpha\gamma-\beta^2, \, \beta\delta-\gamma^2, \, \alpha\delta-\beta\gamma) \eeq
and hence $V_A\cong \P^1$ again, but embedded into $\P^3$ as a twisted cubic curve:
$$ x\!:\! y \; \mapsto \; x^3\!:\! x^2y \!:\! x y^2 \!:\!y^3 $$
The generalisation to higher degree polynomials $f(X)$ is obvious. 

\end{egu}

Note that the inclusion
$$i: V_A \into \P^{n-1}$$
 is a \emph{toric embedding}, the image is not invariant under the full torus action on $\P^{n-1}$ but it is preserved by a subtorus of rank $d$. Also the moment polytope of $V_A$ is just the convex hull of $A$ (see Section \ref{fromchow}). 
\pgap

Now choose a hyperplane in $H\subset \P^{n-1}$. It is specified by a section of $O(1)$, a linear function with some coefficients $a_1,..., a_n$. Pulling this section back along $i$ we get a section of a line bundle of $V_A$, and in the big open torus $(\C^*)^d \subset V_A$ this is exactly our Laurent polynomial $f_A$.  Obviously 
$$H \mbox{ is tangent to }V_A \; \iff\; (f_A) \mbox{ has a singular point} $$ 
and we conclude that $\Delta_A$ is the equation for the projective dual of $V_A$. 

\begin{rem} In this argument we're ignoring points of $V_A$ that live outside $(\C^*)^d$. This makes no difference when computing the projective dual; the tangent hyperplanes to points in $V_A\cap (\C^*)^d$ provide a Zariski-open subset of $(V_A)^\vee$, and we know that $(V_A)^\vee$ is always irreducible. \end{rem}

\begin{eg} \label{quadcurve2} Set $A=\{X^2, X, 1\}$ so $V_A\cong \P^1$ is a quadric curve in $\P^2$ (Example \ref{quadcurve1}).  A hypersurface $H$ is defined by a linear form $a\alpha+b\beta+c\gamma$, and pulling this back to $V_A$ gives a section of $O(2)$, the homogeneous quadratic form:
$$a X^2 + bXY + cY^2 $$
In the torus $\C^*\subset V_A$ this is $f_A(X)$. So the projective dual to $V_A$ is:
$$V_A^\vee = (b^2-4ac)\subset \P^2$$
 Hence the projective dual to any quadric curve is another quadric curve.
\end{eg}

\begin{eg} \label{twistedcubic2} The projective dual to a twisted cubic $\P^1\subset \P^3$ (Example \ref{twistedcubic}) is a quartic surface in $\P^3$, defined by the discriminant $\Delta_3$ of a cubic polynomial \eqref{cubicdis}.

 Note that this must be a special quartic surface since its projective dual is a curve. If we take the projective dual of a general quartic surface (or even a general quartic having the same Newton polytyope as $\Delta_3$) we will get a surface in $\P^3$ of some high degree.
\end{eg}
\pgap

\subsection{Toric degenerations}

Finding projective duals or discriminants directly is hard in most cases. An easier question is to look only for the extremal terms of the discriminant, as discussed in Section \ref{sec.intro}. This means using the torus action on projective space and studying projective duality `asymptotically'. 
\pgap

Let $V\subset \P^n$ be any subvariety of projective space. If we act by the torus $(\C^*)^n$ on $\P^n$ we generate an $n$-parameter family of subvarieties, and the limits of this family (once we've defined them properly) must be subvarieties which are torus-invariant. We will call these limits the \emph{extreme toric degenerations} of $V$. 
\pgap

The straightforward case is when $V$ is a hypersurface.

\begin{egu} Let $V\subset \P^2_{\alpha:\beta:\gamma}$ be the quadric curve cut out by:
$$F(\alpha, \beta, \gamma) = \alpha^2 + \alpha\beta + \beta^2 + \alpha\gamma + \beta\gamma$$
The torus action on $\P^2$ rescales the coefficients of $F$, moving them away from 1,  producing a 2-parameter family of quadric curves. It's easy to see that the limits of this family are the four reducible and torus-invariant quadrics
$$(\alpha^2),\quad\quad (\beta^2), \quad\quad (\alpha\gamma),  \quad\quad (\beta\gamma) $$
corresponding to the four vertices of the Newton polytope of $F$. Indeed, a 1-parameter subgroup of $(\C^*)^2$ corresponds to a choice of dual vector $\mathbf{n}$ on the lattice $\Z^2$, and the limit of of $V$ under that $\C^*$ action corresponds to the vertex which maximizes $\mathbf{n}$. So these are the extreme toric degenerations of $V$.

Note that our $V$ in this example is not generic because $F$ has no $\gamma^2$ term, \ie $V$ contains the torus fixed point $0\!:\!0\!:\!1$. A generic quadric curve would have three extreme toric degenerations $(\alpha^2), (\beta^2), (\gamma^2)$. 

\end{egu}

It's clear that all hypersurfaces behave in this way. We can compactify the space of irreducible hypersurfaces using the space of effective divisors in $\P^n$, \ie by  (the projectivization of) the appropriate space of polynomials. And if we have a hypersurface $V=(F)$ then we have a correspondence:
$$\{\mbox{Extreme toric degenerations of $V$}\} \;\leftrightarrow\; \{\mbox{Extremal terms of }F\} $$
When $V$ is not a hypersurface the situation is more complicated. Obviously each extreme toric degeneration must be a union of torus-invariant linear subspaces, since no other subvarieties are torus-invariant. And if $V$ is generic then it's not too hard to argue that every limit must actually be a single linear subspace, perhaps with some multiplicity. But if $V$ is non-generic then taking the limits requires more thought. See Example \ref{twistedcubic3}.
\pgap

Now suppose $V_A\subset \P^n$ is one of our torically-embedded subvarieties from the previous section. This is far from generic, in at least two ways:
\begin{enumerate}\setlength{\itemsep}{5pt}
\item $V_A$ is preserved by a subtorus of rank $d=\dim V_A$ so we only get an $(n-d)$-parameter family of subvarieties. 
\item $V_A$ contains at least $d+1$ of the torus fixed points in $\P^n$, corresponding to the vertices of $A$. Hence all subvarieties in the family, including the limits, must contain all these fixed points. 
\end{enumerate}

\begin{eg} \label{quadcurve3} Let $A=\{X^2,X,1\}$ as in Example \ref{quadcurve1}. We can think of $V_A$ either as the hypersurface $(\alpha\gamma - \beta^2)\subset \P^2$ or as the image of $\P^1$ under the embedding \eqref{eq.quadembed}.

The extreme toric degenerations of $V_A$ are the reducible toric-invariant curves cut out by $\alpha\gamma$ and $\beta^2$. If we take the embedding viewpoint then the second limit corresponds to the double-cover
 $$x\!:\! y \;\mapsto\; y^2\!:\! 0 \!:\! x^2$$
of $(\beta)$. But in the first limit our parametrized curve breaks into two separate components $(\alpha)$ and $(\gamma)$, and the limit is not parametrized by a single $\P^1$. 

Here's a schematic drawing of this family with the original $V_A$ in blue. In one limit we get two sides of the triangle, and in the other we get the remaining side with multiplicity two. 
%\setlength\intextsep{5pt}
%\begin{figure}[h]
\begin{center}
  \includegraphics[width=0.4\textwidth]{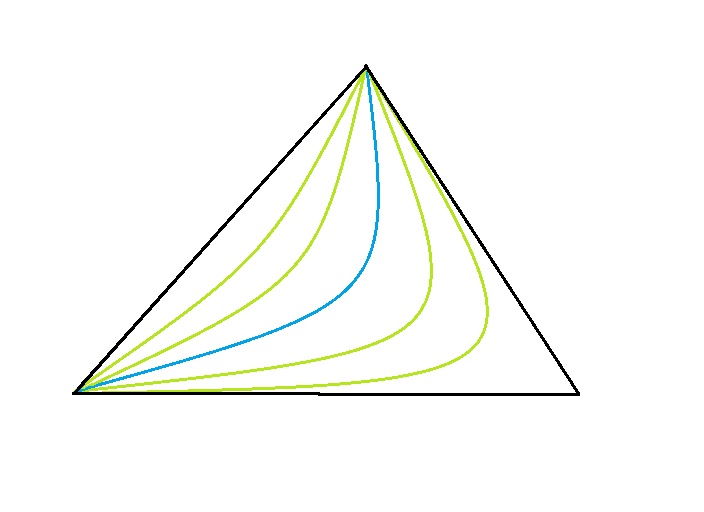}
\end{center}
%\vspace{-20pt}
%\caption{}\label{quadriccurve}
%\end{figure}
\end{eg}

\begin{eg} \label{quadsurf2} Let $A=\{1, X, Y, XY\}$ as in Example \ref{quadsurf}. Then $V_A\cong \P^1\times \P^1$, embedded into $\P^3$ via the Segre embedding
$$ (x\!:\!z, \,y\!:\!w)\; \mapsto\; zw\!:\! xw \!:\! zy \!:\! xy $$
with image cut out by:
$$F(\alpha, \beta, \gamma, \delta) = \alpha\delta - \beta\gamma$$
Notice that $V_A$ contains all four toric fixed points and also two of the toric boundary curves. Schematically we can think of it as a square (the convex hull of $A$) embedded into a tetrahedron, as in this picture:
\begin{center}
  \includegraphics[width=0.4\textwidth]{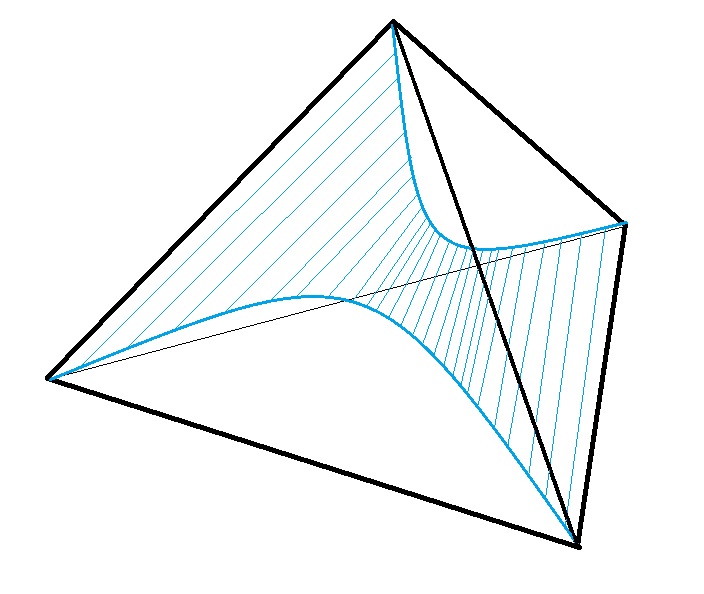}
\end{center}
Acting via $(\C^*)^3$ produces a 1-parameter family of subvarieties in $\P^3$, and the extreme toric degenerations are the reducible toric-invariant surfaces $(\alpha\delta)$  and $(\beta\gamma)$. In our schematic picture the embedded square flows towards the boundary of the tetrahedron, and in the limit breaks into a union of two boundary triangles. The two possible limits correspond to the two possible triangulations of the square. 
\end{eg}

\begin{egu} \label{twistedcubic3} Let $A=\{X^3, X^2, X, 1\}$ as in Example \ref{twistedcubic}, so $V_A\cong\P^1$ embedded as the twisted cubic in $\P^3$. Acting by $(\C^*)^3$ produces a 2-parameter family of embedded curves
$$ x\!:\! y \; \mapsto\; x^3\!:\!  sx^2y \!:\!  txy^2\!:\! y^3 \;\in \P^3$$
each containing the fixed points $1\!:\!0\!:\!0\!:\!0$ and $0\!:\!0\!:\!0\!:\!1$. Since this $V_A$ is not a hypersurface defining the  limits of this family is a bit harder, and we'll wait until the next section to treat this properly. But we can do some heuristic investigation to find that there are four limits:

\begin{itemize}\setlength{\itemsep}{5pt}
\item  In the limit $(s,t)\to (0,0)$ we get a triple cover of the toric boundary curve $\P^1_{\alpha:\delta}\subset\P^3$. 

\item If we send $s\to 0$ we get a 1-parameter family of curves lying in the boundary plane $\P^2_{\alpha:\gamma:\delta}$, defined by $(\gamma^3 - t^3\alpha\delta^2)$. At $t\to 0$ we get $\P^1_{\alpha:\delta}$ with multiplicity three, as above. At $t\to \infty$ we get $\P^1_{\alpha:\gamma}\cup \P^1_{\gamma:\delta}$ where the first component has multiplicity two.
\item Similarly, setting $t\to 0$ gives a 1-parameter family having a limit $\P^1_{\alpha:\beta}\cup \P^1_{\beta:\delta}$ where the second component has multiplicity two.
\item With a little geometric imagination we can see that the final limit is:
$$\P^1_{\alpha:\beta}\cup \P^1_{\beta:\gamma}\cup \P^1_{\gamma:\delta}$$ 
This appears when we send $s\to\infty$ and $t\to\infty$ simultaneously. An algebraic way to see it is to take the ideal of $V_A$ \eqref{twistedcubicideal} and observe that it can degenerate to $(\alpha\gamma, \beta\delta, \alpha\delta)$.\footnote{Although we are actually discussing limits in the Chow variety, as we'll see in the next section, and taking limits of ideals would mean working in the Hilbert scheme. In the Hilbert scheme there are more possible limits but some of them get identified under the Hilbert-Chow morphism.}

\end{itemize}
Schematically the original $V_A$ is an interval embedded into a tetrahedron as a curved path joining two of the vertices. When we flow it to the boundary it breaks into unions of edges. There are four ways this can occur and they correspond to the four subdivisions of the interval we saw in Section \ref{cubic}.

\end{egu}

We can perform this procedure for any of our torically-embedded varieties $V_A$. In the limits $V_A$ will break into unions of parts of the toric boundary of $\P^n$, with some multiplicities.  Given the previous examples, the following claim should be very plausible.

\begin{claim}\label{claim.degenerations} The extreme toric degenerations of $V_A$ correspond to triangulations of the convex hull of $A$, and the multiplicity of each component matches the volume of the corresponding simplex. 
\end{claim}

This is perhaps the key geometric ingredient in Theorem \ref{mainthm}. Indeed we now have enough to attempt a sketch proof of the theorem, as follows:
\begin{align*} \{\mbox{Triangulations of }A\} & \leftrightarrow \{\mbox{Extreme toric degenerations of }V_A\}\\
& \leftrightarrow \{\mbox{Extreme toric degenerations of }V_A^\vee\}\\
& \leftrightarrow \{\mbox{Extremal terms of }\Delta_A\}
\end{align*}
Of course there are still glaring issues to address, including the fact that we haven't defined extreme toric degenerations rigorously yet (beyond the hypersurface case), and the difference between $\Delta_A$ and $E_A$. But there is also important flaw in this sketch in the second step, because taking projective duals does not behave well with respect to toric degenerations. 

\begin{eg} We saw in Example \ref{quadcurve3} that the projective dual of the quadric curve $V=(\beta^2-\alpha\gamma)\subset \P^2$ is another quadric curve $V^\vee=(b^2-4ac)$. Also, one of the extreme toric degenerations of $V$ is the reducible linear curve $V_\infty=(\alpha\gamma)$.  

The projective dual of $V_\infty$ is obviously two points, namely the two hypersurfaces $(\alpha)$ and $(\gamma)$ viewed as points:
$$\{ 1\!:\! 0 \!:\! 0, \; 0\!:\! 0 \!:\! 1\} \subset (\P^2)^\vee$$
 But the extreme toric degenerations of $V^\vee$ are curves. 
\end{eg}

This example illustrates a general phenomenon: an extreme toric degeneration is always a union of linear subspaces, so its projective dual will never be a hypersurface.  We will repair this in the next section by replacing projective duals with another related notion.

\begin{comment}
Of course we haven't yet given a proper definition of an extreme toric degeneration because we haven't identified the correct space in which to take these limits. But in the special case when $V_A$ is a hypersurface there is no issue; the limits exist as divisors in $\P^n$. We saw this in Examples \ref{quadcurve3} and \ref{quadsurf2}.  Moreover, if the defining polynomial of $V_A$ is $g_A(\alpha_1,..., \alpha_{n+1})$ then there is another more obvious bijection:
$$\{\mbox{Extreme toric degenerations of $V_A$}\} \leftrightarrow \{\mbox{Extremal terms of }g_A\} $$
This is a glimpse of another key idea in the proof of Theorem \ref{mainthm} but we'll have to  wait until the next section before we can use it properly. Here it is rather trivial: since $V_A$ is a hypersurface then acting by $(\C^*)^n$ produces only a 1-parameter family, so there are only two extreme toric degenerations. On the other side we can see that $g_A$ will be a polynomial with exactly two terms (having coefficients 1 and $-1$), which encodes the single relation that must hold in the monoid $\Gamma_A$.
\end{comment}

\subsection{Associated hypersurfaces}\label{associated}

Let $V\subset \P^{n-1}$ be a subvariety of projective space of dimension $d$. We've seen how to define the projective dual $V^\vee$ by taking linear hypersurfaces and seeing which ones fail to intersect $V$ generically. There are many variations on this idea. 

For instance, instead of the dual projective space we can consider the Grassmannian $\Gr(n-d-1, n)$ of linear subspaces $\P^{n-d-2}\into \P^{n-1}$ of codimension $d+1$. A generic such subspace will not intersect $V$ at all, but we can consider the subset of those which do intersect $V$, and by another dimension count we see that we expect this subset to be a hypersurface in the Grassmannian. We call it the \emph{associated hypersurface}:
$$ Z_V = \big\{ S \in \Gr(n-d-1, n)\;\; |\; \; S\cap V \neq \emptyset \big\} $$
This construction is also a duality, the variety $V$ can be recovered from its associated hypersurface $Z_V$ \cite[p102]{GKZ}. And it is better behaved than projective duality in that the naive dimension count essentially never fails; if $V$ is irreducible then $Z_V$ will always be an irreducible hypersurface in the Grassmannian \cite[p99]{GKZ}. 

We extend the definition to reducible varieties by treating them as formal sums of irreducible varieties, \ie as effective algebraic cycles, and then the associated hypersurface is the corresponding effective divisor in the Grassmannian.
\pgap

We are interested in the equation cutting out $Z_V$. However, we need to explain exactly what we mean by this, and in fact there are two slightly different options. 
\pgap

A point in $S\in\Gr(n-d-1, n)$ is a subspace defined by $d+1$ linear equations, \ie by a matrix
$$M \in \Mat_{n\times (d+1)}(\C)$$
 having full rank.  The Grassmannian itself is the quotient of this space of matrices by the group $GL_{d+1}(\C)$. This is of course the canonical isomorphism 
$$\Gr_{sub}(n-d-1, \,n)\; \cong \; \Gr_{quot}(n, \,d+1)$$
between a Grassmannian of subspaces of $\C^n$ and a Grassmannian of quotient spaces of $\C^n$. So our associated hypersurface $Z_V$ is the zero locus of some polynomial in $n\times (d+1)$ variables (the entries of $M$) which is invariant under $GL_{d+1}(\C)$. We call this polynomial the $V$\emph{-resultant} and denote it by $\cR_V$.\footnote{\cite{GKZ} use the notation $\tilde{R}_V$. The reason for the name will become clear shortly.}
\pgap

A basic fact about Grassmannians is that any such $\cR_V$ must be expressible as a polynomial in the maximal minors of $M$. These are the Pl\"ucker co-ordinates, \ie we use the Pl\"ucker embedding
$$\Gr(n, d+1)\; \into\; \P(\Wedge^{d+1}\C^n) $$
and then $Z_V$ must in fact be the intersection of the Grassmannian with a hypersurface in the ambient projective space. So a second point-of-view is that $Z_V$ is defined by a polynomial in $n \choose d+1$ variables. This is called the \emph{Chow form}  and we denote it $R_V$.

\begin{eg}\label{eg.linear} Suppose $V$ is a linear subspace $\P^k\subset \P^n$. It's easy to see that $Z_V$ is then a linear hypersurface, \ie $R_V$ is a linear function in the Pl\"ucker co-ordinates. So $\cR_V$ is a polynomial of degree $d+1$. 

Note that individial Pl\"ucker co-ordinates correspond to the linear subspaces $\P^k\subset \P^n$ which are torus invariant.
\end{eg}

\begin{eg} \label{eg.Chowofhyper} Suppose $V$ is a hypersurface in $\P^{n-1}$. Then obviously $Z_V$ is just $V$ itself. So $R_V$ is just the equation defining $V$, under the identification of $\P^{n-1}$ with $\Gr(n, n-1)$. But $\cR_V$ is an equation of much higher degree, obtained by replacing each co-ordinate in $\P^{n-1}$ with the corresponding minor of $M$. 
\end{eg}

Now we specialize this construction to the case when $V$ is one of our torically-embedded subvarieties:
$$i: V_A \into \P^{n-1}$$
Each column of our matrix $M$ is a section of $\cO(1)$ on $\P^{n-1}$, \ie a linear hypersurface, and the intersection of these hypersurfaces is our subspace $S$. Pulling these sections back along $i$ we get $(d+1)$ sections of a line bundle on $V_A$. In the big open torus $(\C^*)^d\subset V_A$ these are Laurent polynomials $g_1, ..., g_{d+1}$, all involving the same set of monomials $A$.  It is evident that
$$ V_A \cap S \neq \emptyset \; \iff \; \mbox{$g_1, ..., g_{d+1}$ have a common root.} $$

So the polynomial $\cR_A:=\cR_{V_A}$ is exactly the \emph{resultant} of this set of polynomials, see Remark \ref{resultants}.

\begin{rem}\,\newline\vspace{-10pt}

\begin{enumerate}\setlength{\itemsep}{5pt} \item Since the polynomials $g_1, ..., g_{d+1}$ all involve the same set of monomials this is a \emph{homogenous} resultant. The general case, involving different sets of monomials, is called a \emph{mixed} resultant.

\item  Previously we weren't exactly clear about where the common root should lie but now we can be precise: the resultant $\cR_A$ vanishes iff $g_1,..., g_{d+1}$ have a common root \emph{in the toric variety} $V_A$. This requires making the Laurent polynomials into homogenous polynomials in the right way, see Section \ref{fromchow}.
\end{enumerate}
\end{rem}

\begin{eg}\label{eg.quadricresultant}
Let $A=\{X^2, X, 1\}$ so $V_A\cong \P^1$ is the quadric curve $(\alpha\gamma-\beta^2)\subset \P^2$. We consider 
$$S \in \Gr(1, 3) \cong \Gr(3,2)$$ 
defined by a matrix:
$$M = \mat{ a_1 & a_2 \\ b_1 & b_2 \\ c_1 & c_2 } $$
The columns of $M$ are two linear forms on $\P^2$, which on $V_A$ become two quadratic forms:
$$g_1 = a_1 X^2 + b_1 XY + c_1Y^2, \quad \quad g_2 = a_2 X^2 + b_2 XY + c_2Y^2 $$
So $S$ intersects $V_A$ iff these two quadratics have a common root in $\P^1$. 

But this is of course a special case of  Example \ref{eg.Chowofhyper}: since $V_A$ is a hypersurface $S$ is just a point, and the associated hypersurface $Z_A := Z_{V_A}$ is just $V_A$ itself. Explicitly we have that 
$$S\;=\;(b_1c_2-c_1b_2): (c_1a_2 - a_1c_2) : (a_1b_2 - b_1a_2) $$
and substituting this into the equation for $V_A$ gives us the resultant:\beq{quadres}\cR_A = \big(b_1c_2-c_1b_2\big) \big(a_1b_2 - b_1a_2\big) - \big(c_1a_2 - a_1c_2\big)^2 \eeq
The two quadratics have a common root iff this expression vanishes.
\end{eg}

\begin{egu}
Let $A=\{X^3,X^2, X,1\}$ so $V_A$ is the twisted cubic curve $\P^1 \into \P^3$ (Example \ref{twistedcubic}). The associated hypersurface $Z_A$ lies in $\Gr(2,4) \subset \operatorname{Mat}_{4\times 2} / GL_2 $, it is the subset where two generic cubics
$$ g_1 = a_1 X^3 + b_1 X^2Y + c_1XY^2 + d_1Y^3, \quad \quad g_2 = a_2 X^3 + b_2 X^2Y + c_2XY^2 + d_2Y^3 $$
have a common root. So the polynomial $\cR_A$ is the resultant of the two cubics.

 However, since $V_A$ is not a hypersurface, computing this resultant is not such an easy problem. Indeed if we could solve this problem then we could immediately deduce the discriminant of a single cubic.

\end{egu}

Chow forms give us a way to take limits in the space of subvarieties in $\P^{n-1}$, and in particular to give a rigorous definition to the extreme toric degenerations considered in the previous section. The action of $(\C^*)^{n-1}$ on $\P^{n-1}$ induces a torus action on any $\Gr(k, n)$, and hence on the space of effective divisors in the Grassmannian.  Given a subvariety $V$, suppose that the associated hypersurface 
$$Z_V \subset \Gr(n-d-1,n)$$
has a toric-invariant limit $(Z_V)_\infty$ in the space of effective divisors. We define the corresponding extreme toric degeneration of $V$ to be the subvariety $V_\infty$ whose associated hypersurface is $(Z_V)_\infty$. 
\pgap

We can also express this in terms of the defining equation for $Z_V$. Our torus action lifts to an action of $(\C^*)^n$ on $\Mat_{n\times (d+1)}(\C)$  just by scaling each row of the matrix. Each monomial in $\cR_V$ has a weight under this action, and we call the set of such weights (or perhaps its convex hull) the \emph{Chow polytope} of $V$.  Note that every monomial contributing to a given maximal minor will have the same weight, so it doesn't matter if we use $\cR_V$ or the Chow form $R_V$ for this definition.
\pgap

Even after passing to $R_V$ it can still be the case that several monomials have the same weight, since we are not using the full torus action on $\P(\Wedge^{n-2-1}\C^n)$. So the Chow polytope is not the Newton polytope of $R_V$, but only a projection of it. Nevertheless we claim that the \emph{vertices} of the Chow polytope must correspond to individual monomials in $R_V$. So we can call these the `extremal terms' of $R_V$. 

To understand this claim, suppose that $V$ degenerates to some torus-invariant $V_\infty$. Then $V_\infty$ must be a union of torus-invariant linear subspaces, perhaps with some multiplicities. For each such subspace the associated hypersurface is cut out by a single Pl\"ucker co-ordinate (Example \ref{eg.linear}). So by definition, the associated hypersurface to $V_\infty$ is the divisor given by the corresponding monomial in Pl\"ucker co-ordinates. This monomial is one of the extremal terms of $R_V$. 

\begin{egu}\label{eg.pentagon}
Let $A=\{1, X, X^2, Y, XY\}$.  Then $V_A$ is a toric surface embedded into $\P^4$, the associated hypersurface $Z_A$ lives in $\Gr(2, 5)$, and computing its defining equation $R_A$ is the same problem as finding the resultant of three polynomials $g_1, g_2, g_3$ where:
$$ g_i = a_i  + b_i X + c_iX^2 + d_iY + e_iXY $$
As we argued in the previous section, extreme toric degenerations of $V_A$ correspond to triangulations of $A$. For example the triangulation
$$\begin{tikzpicture}[scale=1]
\draw (0,0)--(1,0)--(2,0)--(1,1)--(0,1)--(0,0);
\draw (0,1)--(1,0);
\draw (2,0)--(0,1);
\filldraw (0,0) circle (1pt) node[below, yshift=0pt]{$a$};
\filldraw (1,0) circle (1pt)node[below, yshift=0pt]{$b$};
\filldraw (2,0) circle (1pt)node[below,yshift=0pt]{$c$};
\filldraw (0,1) circle (1pt)node[above, yshift=0pt]{$d$};
\filldraw (1,1) circle (1pt)node[above, yshift=0pt]{$e$};
\end{tikzpicture}$$
corresponds to the degeneration:
$$(V_A)_\infty = \P^2_{\alpha:\beta:\delta}\cup \P^2_{\beta:\gamma:\delta}\cup\P^2_{\gamma:\delta:\epsilon}\; \;\;\;\subset \P^4_{\alpha:\beta:\gamma:\delta:\epsilon}$$
The linear subspace $\P^2_{\alpha:\beta:\delta}$ has associated hypersurface cut out by the single Pl\"ucker co-ordinate
$$\pi_{abd} := \det \mat{a_1 & a_2 & a_3 \\ b_1 & b_2 & b_3 \\ d_1 & d_2 & d_3 }$$
and similarly for the other two components. So the associated hypersurface for $(V_A)_\infty$ is the divisor in $\Gr(2,5)$ defined by:
$$\pi_{abd}\pi_{bcd}\pi_{cde}$$
Hence this monomial must be one of the extremal terms of $R_A$. There are four other triangulations of $A$ and they give four more extremal terms:
$$\pi_{abd}\pi_{bde}\pi_{bce}, \quad \pi_{ade}\pi_{abe}\pi_{bce}, \quad \pi_{ade}(\pi_{ace})^2, \quad (\pi_{acd})^2\pi_{cde} $$
\end{egu}

We can now formulate Claim \ref{claim.degenerations} more precisely as:
\begin{thm}\cite[p260]{GKZ} \label{thm.chow}
The extremal terms of $R_A$ biject with triangulations of $A$. For a given triangulation, the corresponding monomial is
$$ \pm \prod (\pi_{A'})^{vol(A')} $$
where the product runs over the simplices $A'$ of the triangulation and $\pi_{A'}$ denotes the corresponding Pl\"ucker co-ordinate.
\end{thm}

Although this result has a clear geometric motiviation the actual proof is entirely algebraic. In Section \ref{cayley} we explain the beginnings of the method they use. 

Note that the fact that the extremal terms are monic (up to sign) is important, but we won't attempt to justify it here. 

\subsection{From Chow forms to discriminants}\label{fromchow}

In this section we fill in the final steps connecting Theorem \ref{thm.chow} to Theorem \ref{mainthm}, and in the process discover how the principal $A$-determinant $E_A$ arises in the story.

\begin{eg} In Example \ref{eg.quadricresultant} we found the formula \eqref{quadres} for the resultant of two general quadratic polynomials $g_1, g_2$, which in our new notation we could write as:
$$\cR_A = \pi_{ab}\pi_{bc} - (\pi_{ac})^2 $$
To get the discriminant of a single quadratic $f = aX^2+bX +c$ we compute the resultant of $g_1=f$ and $g_2 = \partial_X\!f =  2aX + b$, \ie we specialize variables
$$\mat{ a_1 & a_2 \\ b_1 & b_2 \\ c_1 & c_2 }  \leadsto \mat{ a & 0 \\ b & 2a \\ c & b } $$
and hence:
$$\pi_{ab}\mapsto 2a^2, \quad \pi_{bc}\mapsto b^2-2ac, \quad \pi_{ac}\mapsto ab$$ 
Substituting this into $\cR_A$ we do indeed recover the discriminant $\Delta_2$ of a quadratic (times a factor of $a^2$). 

However, we can see that this computation would be tidier if we were to use the \emph{logarithmic derivative}
$$ g_2 = X \partial_X\!f = 2aX^2 + bX $$
instead of the usual $\partial_X\!f$.  Then we'd apply a different specialization
$$\mat{ a_1 & a_2 \\ b_1 & b_2 \\ c_1 & c_2 }  \leadsto \mat{ a & 2a \\ b & b \\ c & 0 } \quad \quad \begin{array}{l}\pi_{ab}\mapsto -ab \\  \pi_{bc}\mapsto -bc \\\pi_{ac}\mapsto 2ac\end{array}$$
and the connection to $\Delta_2$ is more obvious. Notice that this with approach the result is actually $ac\Delta_2$, which is precisely $E_A$ \eqref{eq.quadricEA}.
\end{eg}

In more complicated examples this use of the logarithmic derivative is not only tidier, it's actually essential if we want to get any information.

\begin{egu}\label{eg.pentagon2}Let $A=\{1, X, X^2, Y, XY\}$ as in Example \ref{eg.pentagon}. There we considered three polynomials $g_1, g_2, g_3$ with this Newton polygon and studied their resultant $\cR_A$, finding the five extremal terms. 

Now suppose we want to find the discriminant of  a single polynomial:
$$f_A = a + bX + cX^2  + dY + e XY $$
The obvious thing to do is to set $g_1=f_A$ and $g_2= \partial_X\!f_A$ and $g_3=\partial_Y\!f_A$, so we specialize $\cR_A$ under:
$$\mat{ a_1 & a_2 & a_3\\ b_1 & b_2 & b_3 \\ c_1 & c_2 & c_3 \\ d_1 & d_2 & d_3\\ e_1 & e_2 & e_3}  \leadsto 
\mat{ a & b & d\\ b & 2c & e \\ c & 0 & 0 \\ d & e & 0\\ e & 0 & 0 } 
$$

But then the minors $\pi_{ace}, \pi_{bce}$ and $\pi_{cde}$ will all specialize to zero, and unfortunately all the extremal terms of $\cR_A$ contain one of these three Pl\"ucker co-ordinates (since in any triangulation the edge $ce$ must appear in one of the triangles). So this process tells us nothing at all about $\Delta_A$. 

Geometrically, this specialization is a map
$$\P^4 \to \Gr(2,5)$$
which is not equivariant for the torus actions, so it's not surprising that it doesn't work well with our method.
\end{egu}

Before we continue with the previous example let's recall a little more toric geometry, so that we can understand logarithmic derivatives more clearly. 
\pgap

We have defined the toric variety $V_A$ as the closure of a torus embedded in $\P^{n-1}$, or as Proj of the monoid ring $\C[\Gamma_A]$. However, it can also be constructed as a (GIT) quotient of a vector space by a torus, and this description has some advantages.
\pgap

For each codimension 1 face of the polytope $A$ we consider the primitive normal vector, this gives us a set of vectors $\{v_1,..., v_N\}\subset \Z^d$. These are the rays of the toric fan for $V_A$, which is the normal fan to $A$. Together these vectors define a map $\Z^N \to \Z^d$, which has a kernel $K: \Z^{N-d} \into \Z^N$, expressing the relations between the $v_i$'s. 

The variety $V_A$ is a quotient of the vector space $\C^N$ by the action of the torus $(\C^*)^{N-d}$ having weight matrix $K$.\footnote{This is a GIT quotient, so we need to delete a certain closed subset of $\C^N$ before quotienting.} The co-ordinates $X_1,..., X_N$ on $\C^N$ are natural `projective co-ordinates' on $V_A$.

A section of a line bundle on $V_A$ is given by a polynomial in the $X_i$'s which is semi-invariant under $(\C^*)^{N-d}$, \ie it is quasi-homogeneous in $N-d$ different ways. We have seen that our Laurent polynomial $f_A$ becomes a section of such a line bundle, which appears when we pull-back a generic section of $\cO(1)$ via the embedding $i: V_A \into \P^{n-1}$. So there must be natural way to complete $f_A$ to a multi-quasi-homogenous polynomial in $N$ variables. And there is; you assign each element of $A$ a multi-degree by measuring its lattice distance from each codimension 1 face, and multiply the corresponding monomial accordingly. 

\begin{eg} If $d=1$ then we only have one variable and $A$ is an interval of some length $l$. Hence it has two normal vectors $v_1=1$ and $v_2=-1$ with one relation $v_1+v_2=0$. So $N=2$ and $V_A$ is the quotient of $\C^2$ by $\C^*$ acting with weight $(1,1)$. Hence $V_A\cong \P^1$. We saw this already in Example \ref{twistedcubic}.

 The Laurent polynomial $f_A(X)$ becomes a homogeneous polynomial $f_A(X,Y)$ of degree $l$, which is a section of  the line bundle $\cO(l)$ on $\P^1$. 
\end{eg}

\begin{eg}\label{eg.pentagon3}
Take $A=\{1, X, X^2, Y, XY\}$ as in Examples \ref{eg.pentagon} and \ref{eg.pentagon2}. Then $A$ has four normal vectors $v_1,..., v_4\in \Z^2$ which satisfy two relations:
\begin{align*}\begin{tikzpicture}[scale=1, baseline=0cm]
\draw[-latex] (0,0)--(1,0) node[above, yshift=2pt]{$v_1$};
\draw[-latex] (0,0)--(0,1) node[right, yshift=-5pt]{$v_2$};
\draw[-latex] (0,0)--(0,-1) node[right, yshift=5pt]{$v_3$};
\draw[-latex] (0,0)--(-1,-1) node[above, yshift=5pt]{$v_4$};
\end{tikzpicture} & \hspace{1cm}\begin{array}{c} v_2 + v_3 = 0 \\ v_1+v_2 + v_4=0\end{array} &  K= \mat{0& 1 & 1& 0 \\ 1 & 1& 0 & 1}\end{align*}
Hence the toric surface $V_A$ is a quotient of $\C^4$ by the action of $(\C^*)^2$ with weights $K$.  The Laurent polynomial $f_A(X,Y)$ becomes the homogeneous polynomial
$$f_A(X, Y, Z, W) = aZW^2 + bXZW + c X^2 Z + dYW + eXY $$
where the degrees in $X, Y, Z, W$ correspond to the distances along $v_1,v_2, v_3, v_4$. This polynomial is degree one in $(Y,Z)$ and degree two in $(X, Y, W)$, and it defines a section of a particular line bundle on $V_A$. 
\end{eg}

Once we have homogenized $f_A$ as above we have $N$ different logarithmic derivatives $X_i\partial_{X_i} f_A$  we can consider, and each one has the same quasi-homogeneity properties as $f_A$. So they are sections of the same line bundle that $f_A$ is a section of.\footnote{Geometrically there is no differentiation operation on sections of a line bundle without some additional structure. Here we are using the structure of $V_A$ being a toric variety (or at least the choice of toric boundary).}

Now suppose we look for a common root of $f_A$ and all its logarithmic derivatives. These are $N+1$ polynomials but they are not linearly independent; each quasi-homogeneity of $f_A$ is exactly a relation between them. So we are only really asking for a common root of $d+1$ polynomials, any choice of basis for the lattice:
\beq{alllogd}\big\langle f_A, \;  \;X_1\partial_{X_1}\! f_A, \;\; ..., \;\;  X_N\partial_{X_N}\! f_A \big\rangle_{\Z} \eeq
As usual there should be a hypersurface
$$\nabla \subset (\P^{n-1})^\vee$$
in the space of coefficients such these $d+1$ polynomials have a common root in $V_A$. The equation for $\nabla$ is the resultant of the polynomials.
\pgap

Moreover there are different ways to solve the set of equations \eqref{alllogd} so we expect $\nabla$ to have several irreducible components.

 One option is to find a critical point of $f_A$: this is the projective dual of $V_A$, the vanishing of the discriminant $\Delta_A$. We call this the \emph{principal component} of $\nabla$ (as mentioned in Section \ref{sec.intro}).

But we can also solve the equations by first setting some subset of the $X_i$'s to zero, and then looking for a critical point of the polynomial $f_A$ in the remaining variables. Setting an $X_i$ to zero corresponds to looking at a codimension 1 face $A'\subset A$ of the polytope and considering the associated polynomial $f_{A'}$, or more geometrically, it means restricting $f_A$ to the toric boundary divisor $V_{A'} =\{X_i=0\}\subset V_A$. Then asking for a critical point means looking at the projective dual of $V_{A'}$, \ie the vanishing of the discriminant $\Delta_{A'}$. If we set several $X_i$ to zero we are considering a face $A'\subset A$ of higher codimension.\footnote{We might choose a set of $X_i$ such that the intersection of the corresponding codimension 1 faces of $A$ is empty. But in this case the intersection of the corresponding toric boundary divisors in $V_A$ is also empty, so our polynomials won't have a common root.}

 Putting all the components together, we see that $\nabla$ is exactly the vanishing locus of the principal $A$-determinant $E_A$ \eqref{principalAdet}, so $E_A$ is the resultant of the polynomials \eqref{alllogd}.

\begin{egu}\label{eg.pentagon4}Now we can finish Example \ref{eg.pentagon3}.  We are interested in finding common roots in $V_A$ of the five polynomials:
\beq{alllogd2}\big\{ f_A, \; X\partial_X\! f_A, \; Y\partial_Y \!f_A, \; Z\partial_Z\! f_A, \; W\partial_W \!f_A \big\}\eeq
And we have that 
$$Y\partial_Y\!f_A + Z \partial_Z\!f_A = f_A \aand X\partial_X\!f_A + Y \partial_Y\! f_A + W\partial_W\! f_A   = 2f_A$$
so the vanishing of any three implies the vanishing of all five. There are 6 ways solve these equations, so $\nabla$ has 6 irreducible components. We can either:
\begin{itemize}\setlength{\itemsep}{3pt} \item Find a critical point of $f_A$. This gives the principal component $\{\Delta_A=0\}$. 
\item Set $Y=0$ and find a critical point of $f_{A'} = aW^2 + bXW + cX^2$.
\item Set  $X=Y=0$ and $a=0$. Or the corresponding solution for the other three vertices of A. 
\end{itemize}
So $\nabla$ is cut out by the reducible polynomial:
$$E_A = acde(b^2-4ac)\Delta_A(a,b,c,d,e) $$
On the other hand, we know that $E_A$ can be computed as the resultant of three polynomials, so we can compute its extremal terms using Theorem \ref{thm.chow}. Indeed we already computed the extremal terms of $\cR_A$ in Example \ref{eg.pentagon}.

 If choose our three polynomials to be
\beq{g_i}g_1=f_A, \quad g_2= X\partial_X\!f_A, \quad g_3=Y\partial_Y \!f_A\eeq
then we are specializing the variables of $\cR_A$ as follows:
$$\mat{ a_1 & a_2 & a_3\\ b_1 & b_2 & b_3 \\ c_1 & c_2 & c_3 \\ d_1 & d_2 & d_3\\ e_1 & e_2 & e_3}  \leadsto 
\mat{ a & 0 & 0\\ b & b & 0 \\ c & 2c & 0 \\ d & 0 & d\\ e & e & e }$$
Under this specialization each Pl\"ucker co-ordinate maps to the product of the corresponding variables, multiplied by some integer. For example:
$$\pi_{acd} \mapsto acd \times  \det\! \mat{1 & 0 & 0 \\ 1 & 2 & 0 \\ 1 & 0 & 1 } = 2acd $$
Two facts are clear:
 \begin{enumerate}\setlength{\itemsep}{3pt} \item These expressions are independent (up to sign) of our particular choice of polynomials \eqref{g_i},  they only depend on the lattice spanned by \eqref{alllogd2}.
\item The coefficient of each Pl\"ucker co-ordinate is the volume of the corresponding simplex in $A$.
\end{enumerate}
Combining (2) with Theorem \ref{thm.chow}, we see that the extremal terms of $E_A$ are computed by precisely the combinatorial game \eqref{game2} introduced in Section \ref{sec.intro}.

\end{egu}

The last two facts, and hence the conclusion, are evidently true in all examples. This completes our sketch proof of Theorem \ref{mainthm}.

\begin{rem} In the previous example it's actually not hard to compute $\Delta_A$ by hand. It's given by
$$\Delta_A = ae^2 - bde + cd^2 $$
\ie the hypersurface $(f_A)$ is singular iff the root of the linear equation $d + eX$ also satisfies the quadratic $a + bX + cX^2$. So the extremal terms of $E_A$ can be computed directly and Theorem \ref{mainthm} verified in this example. We invite the reader to try this exercise.
\end{rem}

\section{Other topics}

\newcommand\mX{\breve{X}}
\newcommand\mY{\breve{Y}}
\newcommand\mW{\breve{W}}

In this final section we cover a couple of additional topics that may be of interest. In the previous sections we tried to assume as little background as possible; from now on we make no such attempt. 

\subsection{The Cayley method}\label{cayley}

In 1848 Cayley discovered a beautiful method for computing resultants exactly. It is not computationally efficient but it provides an essential theoretical tool for the proofs in \cite{GKZ}. 
\pgap

Let $E$ be a vector space of dimension $n$ and $s\in E$ an element. The \emph{Koszul complex} associated to $s$ is the chain complex:
$$ \Wedge^n E^\vee \stackrel{s}{\To} ... \stackrel{s}{\To} \Wedge^2 E^\vee \stackrel{s}{\To} E^\vee \stackrel{s}{\To} \C $$
This complex is exact, unless $s=0$. More generally if $E$ is a vector bundle over some base variety $B$, and $s$ is a section of $E$, then the Koszul complex $(\Wedge^\bullet E^\vee, s)$ is a chain complex of vector bundles over $B$. It is exact away from the vanishing locus of $s$. In fact if $s$ is a transverse section then the Koszul complex has homology only at the final term, so it is a free resolution of the torsion sheaf $\cO_{s=0}$ on $B$. 

\begin{eg}
Let 
$$B=\C^6_{a_1,b_1,c_1, a_2, b_2, c_2}\times \P^1_{x:y}$$
and set $E=\cO(2)^{\oplus 2}$. We have a tautological section of $E$ given by two generic quadratics:
$$s = (g_1, g_2)= (a_1 x^2 + b_1xy + c_1 y^2, \, a_2 x^2 + b_2xy + c_2 y^2) \; \in \Gamma (B,  E) $$
The vanishing locus $\{s=0\}$ is the set of points where $x\!:\!y$ is a root of both $g_1$ and $g_2$. The Koszul complex $(\Wedge^\bullet E^\vee, f)$ is a free resolution of the sheaf $\cO_{s=0}$. 
\pgap

Now consider the projection $\pi: B \to \C^6$. The push-down $\pi_*\cO_{s=0}$ is a torsion sheaf, supported on the locus in parameter space where $f_1$ and $f_2$ have a common root. So if we write $\cR_2$ for the resultant of the two quadratics then this sheaf is supported on the subvariety:
$$\{\cR_2 = 0\} \subset \C^6$$
 To get a free resolution of this sheaf we can take the derived push-down of the Koszul complex on $B$, \ie take the cohomology of each bundle $\Wedge^p E^\vee$ over the $\P^1$ fibres. We can make life easier by first twisting by the line bundle $\cO(3)$, so the Koszul complex becomes
$$ \cO(-1) \stackrel{s}{\To} \cO(1)^{\oplus 2} \stackrel{s}{\To} \cO(3) $$
which is a free resolution of the sheaf $\cO_{s=0}(3)$. Pushing down this sheaf gives, again, some torsion sheaf supported on $\{\cR_2=0\}$, and it has a free resolution
\beq{D} 0 \To \cO^{\oplus 4} \stackrel{D}{\To}  \cO^{\oplus 4} \To \pi_*(\cO_{f=0}(3)) \eeq
because:
\begin{itemize}\setlength{\itemsep}{5pt}
\item $R\pi_*\cO(3) = H^0(\P^1, \cO(3)) \otimes \cO_{\C^6} \cong (\cO_{\C^6})^{\oplus 4}  $
\item $R\pi_*\cO(1) = H^0(\P^1, \cO(1)) \otimes \cO_{\C^6} \cong (\cO_{\C^6})^{\oplus 2}  $
\item $R\pi_*\cO(-1) = 0$, since $H^p(\P^1, \cO(-1))=0$ for all $p$. 
\end{itemize}
The differential $D$ in \eqref{D} is some $4\times 4$ matrix of polynomials in $a_1,..., c_2$. And since the cokernel of $D$ is supported on the locus where $D$ drops in rank, we conclude:
$$\det D = \cR_2$$
This provides another way to compute the resultant $\cR_2$. Indeed, using the obvious bases $X,Y$ for $H^0(\P^1, \cO(1))$ and $X^3, X^2Y, XY^2, Y^3$ for $H^0(\P^3, \cO(3))$ we can easily compute that 
$$D = \mat{ 
a_1 & b_1 & c_1 & 0 \\
0 & a_1 & b_1 & c_1  \\
a_2 & b_2 & c_2 & 0 \\
0 & a_2 & b_2 & c_2  }$$
and it's easy to verify that $\det D$ agrees with the expression \eqref{quadres}.
\end{eg}

We can use the same technique to compute resultants of higher degree.

\begin{egu} 
Let $\cR_d$ denote the resultant of a pair of degree $d$ polynomials:
$$g_1(X,Y), \quad g_2(X,Y)$$
So $\cR_d$ is a polynomial in $2(d+1)$ variables.  We set
$$B = \C^{2(d+1)}\times \P^1$$
and then $s=(g_1,g_2)$ is a section of a vector bundle $E=\cO(d)^{\oplus 2}$. We take the Koszul complex of $s$, twisted by $\cO(2d-1)$, to get an exact sequence
$$0 \To \cO(-1) \stackrel{s}{\To} \cO(d-1)^{\oplus 2} \stackrel{s}{\To} \cO(2d-1) \To \cO_{s=0}(2d-1) \To 0 $$
and then pushing-down to the coefficient space $\C^{2(d+1)}$ gives a free resolution:
$$ 0\To \cO^{\oplus 2d} \stackrel{D}{\To} \cO^{\oplus 2d} \To \pi_*(\cO_{s=0}(2d-1))  \To 0 $$
Here $\pi_*(\cO_{s=0}(2d-1))$ is some torsion sheaf supported on $\{\cR_d=0\}$. We conclude that $\det D$ must equal $\cR_d$ (or at least some power of $\cR_d$). 

 For example, setting $d=3$ we get the resultant of two cubics as:
$$\cR_3 = \det \mat{ 
a_1 & b_1 & c_1 & d_1 & 0 & 0   \\
0 & a_1 & b_1 & c_1 & d_1 & 0  \\
0 & 0 & a_1 & b_1 & c_1 & d_1  \\
a_2 & b_2 & c_2 & d_2 & 0 & 0 \\
0 & a_2 & b_2 & c_2 & d_2& 0 \\
0 & 0 & a_2 & b_2 & c_2 & d_2}
$$

By specializing the variables appropriately (as in Section \ref{fromchow}) you can compute the discrimant $\Delta_3$ by hand. For larger $d$ it's an easy problem for a computer algebra system. 
\end{egu}

We can compute more general homogenous resultants $\cR_A$ by the same technique. Take a set of polynomials $g_0,..., g_d$ each using the same set of monomials $A \subset \Z^d$. Then we have a toric variety $V_A$, and after homogenenising correctly each $g_i$ is a section of some line bundle $L$ on $V_A$. We take the space
$$B = \C^{(d+1)|A|} \times V_A $$ 
and consider the tautological section:
$$s = (g_0, ..., g_d) \;\in \Gamma(B, L^{\oplus (d+1)}) $$
The push-down of the Koszul resolution of $\cO_{s=0}$ gives a complex of vector bundles on coefficient space, with homology supported on $\{\cR_A=0\}$. This complex will usually have length $>2$, so we cannot compute $\cR_A$ as the determinant of a single matrix, but it can be computed as the `determinant of a complex' \cite[Appendix A]{GKZ}.

\subsection{Mirror symmetry}\label{mirror}

A modern context in which an algebraic geometer might encounter discriminants is mirror symmetry.\footnote{\emph{E.g.} the author.} This is a huge topic and we'll say very little about it, just enough to explain how discriminants are relevant. 
\pgap

Mirror symmetry is, at heart, an abstract duality acting on a certain class of quantum field theories. In mathematics it is usually a duality on \emph{Landau-Ginzburg models}, this means a pair $(Y,W)$ consisting of a K\"ahler manifold $Y$ and a holomorphic function $W: Y \to \C$. What makes it particularly interesting is that it swaps the roles of the complex and symplectic structures; properties of the complex/algebraic geometry of $(Y,W)$ are reflected by the symplectic geometry of its mirror $(\mY, \mW)$. 

The simplest class of examples involve toric geometry. Suppose we choose $Y$ to be simply a torus
$$Y = (\C^*)^{d}$$
and $W$ some Laurent polynomial. Then $W$ is determined by some set of monomials $A\subset \Z^{d}$, plus the choice of co-efficients. The mirror to $(Y, W)$ will be a toric variety $\mY$ determined combinatorially by $A$, equipped with the zero superpotential $\mW\equiv 0$.\footnote{This is a recipe due originally to Witten and Hori-Vafa. Note that we are not specifying the K\"ahler metric on $\mY$, which is a much harder question.}
\pgap

To build this toric variety $\mY$ we view $A$ as a map $\Z^n \to \Z^{d}$, where $n=|A|$, and find its kernel:
$$Q: \Z^{n-d} \to \Z^n$$
Then we take a GIT quotient of $\C^n$ by the action of the torus $(\C^*)^{n-d}$ having weight matrix $Q$. This gives a $d$-dimensional toric variety whose fan `refines $A$', in the sense that the one-dimensional cones of the fan are the rays spanned by (some subset of) the elements of $A$. 

\begin{rem} This construction of $\mY$  is similar to the construction of $V_A$ described in Section \ref{fromchow}, but it is important to understand that they are different toric varieties. The fan for $V_A$ is the normal fan to the convex hull of $A$, so its rays come from the normal vectors to this polytope, not the elements of $A$.  Also $V_A$ is always compact but $\mY$ may not be. 
\end{rem}

In fact from this data $Q$ there usually be more than one possible GIT quotient we can form because GIT requires a choice of some extra data, a \emph{stability condition}.\footnote{When we build $V_A$ as a GIT quotient there is a canonical stability condition that we use implicity.} So there will be some number of different toric varieties
$$\mY_1,..., \mY_t$$
whose fans are different possible refinements of $A$. These varieties are all birational. 

This non-uniqueness is an important part of the mirror symmetry story. Varying the coefficients of $W$ should correspond to varying the K\"ahler class on the mirror variety, since mirror symmetry exchanges complex moduli for symplectic moduli. And in the space of quantum field theories, large variations of the K\"ahler class can mean passing to another birational model. This leads to the concept of the \emph{stringy K\"ahler moduli space} which unifies the K\"ahler moduli spaces of all the $\mY_i$'s. 

Much of mirror symmetry involves understanding how certain quantities (Gromov-Witten invariants, derived categories, ...) vary over this moduli space. Unfortunately the space itself has no rigorous definition in general, so if we want to compute examples we usually stick to toric varieties, where the stringy K\"ahler moduli space can be \emph{defined} to be the coefficient space of the mirror $W$. So to do these kinds of computations it is essential to understand the discriminant locus of $W$. 
\pgap

As a final remark, this story works particularly cleanly in the case that each $\mY_i$ is a Calabi-Yau variety. This happens when the torus $(\C^*)^{n-d}$ acts inside the the special linear group $SL_n(\C)$, which is a condition on the matrix $Q$, and is equivalent to the monomials $A\subset \Z^{d}$ all living in some affine hypersurface of height 1. This says that $W$ can be written in the form:
$$W = X_0f_A(X_1,...., X_d) $$
In this situation the fans for the $\mY_i$'s correspond exactly to triangulations of the convex hull of $A$ (to turn a triangulation into a fan we just take the cone on each simplex).  So Theorem \ref{mainthm} gives us some relationship between the asymptotics of $E_A$ and the different birational models $\mY_i$ of our toric Calabi-Yau. 

The precise statement involves yet another toric variety associated to $A$. Back in Section \ref{cubic} we observed that the discriminant of a cubic $\Delta_3$ can be sensibly understood as a rational function in two variables, by scaling out some symmetries. By the same process a general $\Delta_A$ can be reduced to a function of $n-d$ variables, by quotienting by a torus of rank $d$. More precisely, we should understand the discriminant locus $\{E_A=0\}$ as a subvariety in some projective toric variety of dimension $n-d$. This is the \emph{secondary toric variety}; its toric data encodes the triangulations of $A$, and in particular it has one toric fixed point for each $\mY_i$. Theorem \ref{mainthm} says that the moment polytope of this secondary toric variety is exactly the Newton polytope of $E_A$, \ie the asymptotics of $\{E_A=0\}$ are dual to the toric boundary.

% ----------------------------------------------------------------
\bibliographystyle{halphanum}

\end{document}